\documentclass{amsart}

\usepackage{fullpage}
\usepackage{graphicx}
\usepackage[tight,footnotesize,FIGTOPCAP]{subfigure}
\usepackage{amsmath}
\usepackage{amsfonts}
\usepackage{amssymb}
\usepackage{amsthm}
\usepackage{bm}
\usepackage{xspace}
\usepackage{psfrag}
\usepackage{color}
\usepackage{ifthen}
\usepackage{upgreek}
\usepackage[nointegrals]{wasysym}
\usepackage{pdflscape}

\theoremstyle{definition}

\theoremstyle{plain}

\newcommand{\setR}{\ensuremath{\mathbb{R}}}

\newcommand{\setC}{\ensuremath{\mathbb{C}}}

\newcommand{\setH}{\ensuremath{\mathbb{H}}}
\newcommand{\setHpu}{\ensuremath{\mathbb{H}_\mathrm{pu}}} 

\newcommand{\mInv}{\ensuremath{^{-1}}}      
\newcommand{\mTr}{\ensuremath{^\mathrm{T}}}      
\newcommand{\mStar}{\ensuremath{^\ast}}          
\newcommand{\mHerm}{\ensuremath{^\mathrm{H}}}    

\renewcommand{\vec}[1]{\ensuremath{\bm{{#1}}}}
\newcommand{\vecTr}[1]{\ensuremath{\bm{#1}\mTr}}
\newcommand{\vecHerm}[1]{\ensuremath{\bm{#1}\mHerm}}

\newcommand{\mat}[1]{\vec{#1}}
\newcommand{\matTr}[1]{\vecTr{#1}}
\newcommand{\matHerm}[1]{\vecHerm{#1}}

\newcommand{\mVec}[1]{\ensuremath{\mathrm{vec}\left( #1 \right)}}    

\newcommand{\mDiag}[1]{\ensuremath{\mathrm{diag}\left( #1 \right)}}  

\newcommand{\mInvL}{\ensuremath{^{\vartriangleleft}}}      
\newcommand{\mInvR}{\ensuremath{^{\vartriangleright}}}      

\newcommand{\mMulL}{\ensuremath{\cdot_\mathrm{\scriptscriptstyle L}}}
\newcommand{\mMulR}{\ensuremath{\cdot_\mathrm{\scriptscriptstyle  R}}}
\newcommand{\mMulLR}{\ensuremath{\cdot_\mathrm{\scriptscriptstyle L/R}}}

\newcommand{\mKronL}{\ensuremath{\otimes_\mathrm{\scriptscriptstyle L}}}
\newcommand{\mKronR}{\ensuremath{\otimes_\mathrm{\scriptscriptstyle R}}}

\newcommand{\mKhatriL}{\ensuremath{\diamond_\mathrm{\scriptscriptstyle L}}}
\newcommand{\mKhatriR}{\ensuremath{\diamond_\mathrm{\scriptscriptstyle R}}}

\newcommand{\hReal}[1]{\ensuremath{\mathrm{s}\left( #1 \right)}}    
\newcommand{\hImag}[1]{\ensuremath{\mathrm{v}\left( #1 \right)}}    

\newcommand{\dontprint}[1]{}

\newcommand{\gleichung}[1]{(\ref{eq:#1})\xspace}

\newcommand{\abschnitt}[1]{Section~\ref{sec:#1}}

\begin{document}


\title{Using Quaternion-Valued Linear Algebra}

\address{Electronic Measurement Research Laboratory, Ilmenau University of Technology, Germany}
\author{Dominik Schulz}
\email{dominik.schulz@tu-ilmenau.de}
\author{Reiner S. Thom\"a}
\email{reiner.thomae@tu-ilmenau.de}


\begin{abstract}
Linear algebra is usually defined over a field such as the reals or complex numbers.
It is possible to extend this to skew fields such as the quaternions. However,
to the authors' knowledge there is no commonly accepted notation of linear
algebra over skew fields.
To this end, we discuss ways of notation that account for the non-commutativity
of the quaternion multiplication.
\end{abstract}

\maketitle

%

\section{Introduction}
\label{sec:introduction}

The use of quaternion linear algebra is emerging among researchers. However, defining
the properties of quaternionic matrices is still subject to research. Nevertheless,
advances have been made especially in analyzing eigenvalues 
(\cite{Zhang97}, \cite{Faren03}, \cite{Zou12}, \cite{Huang01}, \cite{Baker99})
of quaternion matrices. 

Albeit the analytical analyses, a practical framework is needed to use quaternion
matrices. A basic problem arises due to the non-commutativity of the quaternion
algebra. That is, for two matrices $\mat A \in \setH^{M \times K}$ and 
$\mat B \in \setH^{K \times N}$ the matrix products
$\sum_{k=1}^K \left[\mat A\right]_{m, k} \left[\mat B\right]_{k, n}$
and
$\sum_{k=1}^K \left[\mat B\right]_{k, n} \left[\mat A\right]_{m, k}$
are in general not the same.
To this end, we propose a new notation for quaternion matrix multiplication. In this
notation we take care of the multiplication oder. 

Matrices are denoted using capital letters in bold face. Vectors are denoted
by lower case characters in bold face. The $\mTr$ operator indicates the transpose
of a matrix or a vector. Additionally, $\mStar$ and $\mHerm$ stand for the 
conjugation and the Hermitian transpose (conjugate transpose) respectively.
The \mDiag{\cdot} operator transforms a vector into a square matrix having
the elements of the vectors on the main diagonal.
\section{Review on Quaternions}

The current section will introduce the quaternions and the notation used
in the subsequent sections.

%
%
%

\subsection{Quaternions and Complex Numbers}

Quaternions have been first discovered by W. R. Hamilton \cite{Ha1844}.
The quaternions are one of several possible extensions of complex numbers.

The set of quaternions $\setH$ can be constructed from the set of complex numbers $\setC$.
For that, let $z_1 = a_1 + b_1\imath$ and $z_2 = a_2 + b_2\imath$ be two complex numbers. 
Additionally, we introduce another imaginary unit $\jmath$ with $\jmath^2 = -1$.
Then, a quaternion $q \in \setH$ is a pair of complex numbers
$(z_1, z_2) \equiv z_1 + z_2\jmath$.
Hence, we have
\begin{align}
  \label{eq:complex2quaternion}
  q = z_1 + z_2\jmath &= a_1 + b_1\imath + a_1\jmath + b_2 \imath\jmath
\end{align}
By introducing a third imaginary unit $k := \imath \jmath$, with $k^2 = -1$, we obtain the 
common form of a quaternion.
\begin{align*}
  q = a_1 + b_1\imath + a_1\jmath + b_2 k
\end{align*}

From now on, we will denote the set of quaternions as $\setH$. The associated 
algebra is expressed in terms of pairs of complex numbers.
To this end, let $x_1$, $x_2$, $z_1$, and $z_2$ be two pairs of complex numbers that form
two quaternions $q_1 = x_1 + x_2\jmath$ and $q_2 = z_1 + z_2\jmath$.
In this case, addition, multiplication and conjugation are defined as follows.
\begin{align*}
  (x_1, x_2) + (z_1, z_2) &= (x_1 + z_1, x_2 + z_2) \\
  (x_1, x_2) \cdot (z_1, z_2) &= (x_1 z_1 - x_2 z_2^\ast, x_2 z_1^\ast + z_2 x_1) \\
  (x_1, x_2)^\ast &= (x_1^\ast, -x_2)
\end{align*}

\subsection{Symplectic Decomposition}

Let $q = a_0 + a_1\imath + a_2\jmath + a_3 k \in \setH$ be a quaternion. The
value $\hReal{q} = a_0$ is called the \emph{scalar part} (or real part) and the value 
$\hImag{q} = a_1\imath + a_2\jmath + a_3 k$ is called the \emph{vector part}.

A \emph{pure quaternion} is a quaternion having a vanishing scalar part.
Similar to the complex numbers, the modulus of a quaternion $q$ can be expressed 
as $\lvert q \rvert := \sqrt{q^\ast q} = \sqrt{a_0^2 + a_1^2 + a_2^2 + a_3^2}$.

If $q'$ is an arbitrary quaternion with non-vanishing vector part, then
\begin{align*}
  \mu(q') = \frac{\mathrm{v}(q')}{\lvert q' \rvert}
\end{align*}
is a \emph{pure unit} quaternion (PUQ). 
Throughout the subsequent sections we will use $\setHpu$ to indicate the set of
all pure unit quaternions.
\begin{align*}
  \setHpu := \left \lbrace a_1 \imath + a_2 \jmath + a_3 k : \sum_{n=1}^3 a_n^2 = 1,  a_n \in \setR\right \rbrace
\end{align*}
It can easily be proven that for any PUQ $\mu \in \setHpu$ the relation $\mu^2 = -1$ holds. 
Based on this observation, it is possible to build subsets of quaternions that each are isomorphic to complex
numbers. To this end, let $\mu \in \setHpu$ be an arbitrary but fixed pure unit quaternion. Consequently,
the set
\begin{align*}
  \setC_\mu := \left \lbrace a_0  + a_1 \mu : a_n \in \setR\right \rbrace
\end{align*}
is isomorphic to the set of complex numbers $\setC$. One special consequence is
the following fact: Iff for two quaternions $q \in \setH$ and $p \in \setH$ there 
exists a PUQ $\mu \in \setHpu$ such that 
$q \in \setC_\mu$ and $p \in \setC_\mu$, the product of both becomes commutative
(i.e., $q p = p q$).

Before continuing, let us first define the notion of orthogonality of PUQs.
Let $\mu = a_1 \imath + a_2 \jmath + a_3 k \in \setHpu$ and $\mu_\perp = b_1 \imath + b_2 \jmath + b_3 k \in \setHpu$
be two pure unit quaternions. $\mu$ and $\mu_\perp$ are said to be orthogonal iff
$\sum_{n=1}^3 a_n b_n = 0$. 

For notational convenience, in all following analyses we will assume that the Greek letter 
$\mu$ always refers to a PUQ. Additionally, $\mu_\perp$ will always refer to some PUQ which 
is orthogonal to $\mu$. 

Finally, the \emph{symplectic decomposition} as presented in \cite{Ell07} can be obtained 
in the following way: Given a quaternion $q \in \setH$ choose two PUQs $\mu$ and $\mu_\perp$.
Then, $q$ can be decomposed into the following form:
\begin{align}
  \label{eq:PUQ2quaternion}
  q &= q_0 + q_1\mu + (q_2 + q_3\mu)\mu_\perp, & q_n &\in \setR \nonumber \\
   &= q'_0 + q'_1\mu_\perp, & \quad q'_n &\in \setC_\mu 
\end{align}
It turns out that \gleichung{PUQ2quaternion} is as generalized version of \gleichung{complex2quaternion}.
In order to remain most general all subsequent equations will be given in terms
of pure unit quaternions rather than in terms of $\imath$, $\jmath$, and $k$.


\subsection{Euler's Formula}

In the four-dimensional quaternionic space each pure unit quaternion can be seen as an 
axis. In conjunction with the real axis infinitely many planes exist that contain
complex-isomorphic numbers. A rotation within one of such planes is similar to a rotation
in the ordinary complex plane. Hence, Euler's formula obtains the following form
in the quaternion domain.

  Let $q \in \setH$ be some non-zero quaternion. Then, up to a sign ambiguity 
  there exists a unique axis $\mu \in \setHpu$ and an angle $\alpha \in (-\pi, \pi]$ such that
  \begin{align}\label{eq:euler}
    q &= \lvert q \rvert \exp(\mu \alpha) 
      = \lvert q \rvert \left( \cos(\alpha) + \mu \sin(\alpha) \right) .
  \end{align}  
In \gleichung{euler} $\exp(\cdot)$ denotes the exponential function. 
It may be noted that \gleichung{euler} 
also implies the fact that for any quaternion $q$ there exists a $\mu$ such that $q \in \setC_\mu$.

In general, for two quaternions $q$ and $p$ the value $\exp(q+p)$ does not equal $\exp(q)\exp(p)$.
However, if $q$ and $p$ commute, then it is true that $\exp(q+p) = \exp(q)\exp(p)$.

\subsection{Similar Quaternions}

Due to the lack of commutativity in the quaternion algebra the term $p = s^{-1} q s$ for two quaternions
$q, s \in \setH$ is in general not the same as $q$. Nevertheless, $p$ is said to be \emph{similar} to $q$.
This is an equivalence relation, where the equivalence class $\mathcal{S}\left( q \right)$ is defined as
follows.
\begin{align*}
  \mathcal{S}\left( q \right) := 
    \left \lbrace
      p \in \setH
      \colon
      p = s^{-1} q s, s \in \setH, \lvert s \rvert = 1
    \right \rbrace
\end{align*}
\section{Matrix Products}
\label{sec:left_right_mult}

\subsection{Left and Right Matrix Multiplication}

The quaternion multiplication is not commutative. That is, in general $a b$ is not
the same as $b a$ for two quaternions $a \in \setH$ and $b \in \setH$. 
A problem arises when computing the product of two quaternionic matrices 
$\mat A \in \setH^{M \times K}$  and $\mat B \in \setH^{K \times N}$.
\begin{align*}
  \left[\mat A  \mat B\right]_{m, n} &:= \sum_{k=1}^K \left[\mat A\right]_{m, k} \left[\mat B\right]_{k, n}
\end{align*}
The entries of $\mat A$ are multiplied from the left. Likewise, the entries of $\mat B$ are multiplied from the right.
If the order needs to be swapped, the following expression gives the correct result.
\begin{align}
  \label{eq:mat_mul_reverse}
  \left[ \left( \matTr B  \matTr A \right)\mTr \right]_{m, n} &:= \sum_{k=1}^K \left[\mat B\right]_{k, n} \left[\mat A\right]_{m, k}
\end{align}

However, the left-hand side of \gleichung{mat_mul_reverse} is somewhat less intuitive. Hence, in \cite{Schu11} we introduced a notation 
for left multiplication and right multiplication.
\begin{align*}
  \left[\mat A \mMulL \mat B\right]_{m,n} &:= \sum_{k=1}^K \left[\mat A\right]_{m, k} \left[\mat B\right]_{k, n} \\
  \left[\mat A \mMulR \mat B\right]_{m,n} &:= \sum_{k=1}^K \left[\mat B\right]_{k, n} \left[\mat A\right]_{m, k}
\end{align*}
Consider the case of conformant matrices, i.e. where $M=N$. It follows that four different matrix products,
$\mat A \mMulLR \mat B$ and $\mat B \mMulLR \mat A$,
are possible. This highlights the fact that we deal with two different kinds of ordering: The first kind of 
ordering denotes inner products of either the rows of $\mat A$ and the columns of $\mat B$ or inner products
of the rows of $\mat B$ and the columns of $\mat A$. Commonly, this is represented by the the ordering in which 
the matrices $\mat A$ and $\mat B$ appear in an equation. The second kind of ordering refers 
to the ordering of the scalars within each inner product. This ordering is specified by the proposed
operators $\mMulL$ and $\mMulR$.

Next, we observe that the usage of the left and right multiplication operators allows for a convenient description of how the transpose,
conjugation and Hermitian transpose (conjugate transposition) act on a product of two quaternion matrices.
\begin{align}
  \label{eq:qmatTransp}
  \left( \mat A \mMulL \mat B \right)\mTr &= \matTr B \mMulR \matTr A \\
  \label{eq:qmatConj}
  \left( \mat A \mMulL \mat B \right)^\ast &= \mat A^\ast \mMulR \mat B^\ast 
\end{align}
Note that the second relation follows from the fact that for two quaternion scalars $q$ and $p$ we have
$(pq)^\ast = q^\ast p^\ast$.  Combining \gleichung{qmatTransp} and and \gleichung{qmatConj} yields the
following. 
\begin{align*}
  \left( \mat A \mMulL \mat B \right)\mHerm &= \matHerm B \mMulL \matHerm A \\
  \left( \mat A \mMulR \mat B \right)\mHerm &= \matHerm B \mMulR \matHerm A
\end{align*}
Hence, the Hermitian conjugate behaves the same as in the complex case since the type of 
multiplication (left or right) is not altered.

\subsection{Matrix Product of Three Matrices}

Additionally, let us look at the Product of three matrices. There are six different
ways of ordering the product of the matrix elements which are all covered by the operators
introduced above.
\begin{align*}
  \sum_{k=1}^K \sum_{\ell=1}^L a_{m, k} \ b_{k, \ell} \ c_{\ell, n} &= \left[ \mat A \mMulL \mat B \mMulL \mat C \right]_{m, n} \\
  \sum_{k=1}^K \sum_{\ell=1}^L c_{\ell, n} \ b_{k, \ell} \ a_{m, k} &= \left[ \mat A \mMulR \mat B \mMulR \mat C \right]_{m, n} \\
  %
  \sum_{k=1}^K \sum_{\ell=1}^L a_{m, k} \ c_{\ell, n} \ b_{k, \ell}  &= \left[ \mat A \mMulL \left( \mat B \mMulR \mat C \right) \right]_{m, n} \\  
  \sum_{k=1}^K \sum_{\ell=1}^L c_{\ell, n} \ a_{m, k} \ b_{k, \ell}  &= \left[ \left( \mat A \mMulL  \mat B \right) \mMulR \mat C  \right]_{m, n} \\  
  %
  \sum_{k=1}^K \sum_{\ell=1}^L b_{k, \ell} \ a_{m, k} \ c_{\ell, n}  &= \left[ \left( \mat A \mMulR  \mat B \right) \mMulL \mat C \right]_{m, n} \\  
  \sum_{k=1}^K \sum_{\ell=1}^L b_{k, \ell} \ c_{\ell, n} \ a_{m, k}  &= \left[  \mat A \mMulR \left( \mat B  \mMulL \mat C \right) \right]_{m, n} 
\end{align*}
Note that the first two identities do not contain inner brackets due to the following associativity property.
\begin{align}
  \label{eq:identA7}
  \left( \mat A \mMulL \mat B \right) \mMulL \mat C &= \mat A \mMulL \left( \mat B  \mMulL \mat C \right) \\
  \label{eq:identA8}
  \left( \mat A \mMulR \mat B \right) \mMulR \mat C &= \mat A \mMulR \left( \mat B  \mMulR \mat C \right)
\end{align}  
Similarly, iff $\left[\mat A\right]_{m,n} \in \setC_\mu$ and $\left[\mat B\right]_{u,v} \in \setC_\mu$ are in the same complex-isomorphic set $\setC_\mu$,
then associativity does hold for the remaining equations as well.
\begin{align*}
  \left( \mat A \mMulL \mat B \right) \mMulR \mat C &= \mat A \mMulL \left( \mat B  \mMulR \mat C \right) \\
  \left( \mat A \mMulR \mat B \right) \mMulL \mat C &= \mat A \mMulR \left( \mat B  \mMulL \mat C \right)
\end{align*}

\section{The Fundamental Subspaces}
\label{sec:left_right_subspaces}

In the quaternion domain we may define eight fundamental subspaces using the proposed left and right
matrix multiplication. This is due to the fact
that for the row spaces, column spaces and null spaces we have to define the order of multiplication.
Hence, these are the fundamental subspaces of a matrix $\mat A \in \setH^{M \times N}$.

\subsection*{Left row space ($\mathcal{LR}$) and right row space ($\mathcal{RR}$)}
\begin{align*}
  \mathcal{LR}\left( \mat A \right) &:= 
    \left \lbrace
      \vec y \in \setH^{M \times 1}
      \colon
      \vec y = \mat A \mMulL \vec x, \forall \vec x \in \setH^{N \times 1}
    \right \rbrace
\\    
  \mathcal{RR}\left( \mat A \right) &:= 
    \left \lbrace
      \vec y \in \setH^{M \times 1}
      \colon
      \vec y = \mat A \mMulR \vec x, \forall \vec x \in \setH^{N \times 1}
    \right \rbrace    
\end{align*}

%

\subsection*{Left column space ($\mathcal{LC}$) and right column space ($\mathcal{RC}$)}
\begin{align*}
  \mathcal{LC}\left( \mat A \right) &:= 
    \left \lbrace
      \vec y \in \setH^{N \times 1}
      \colon
      \vecTr y = \vecTr x \mMulL \mat A, \forall \vec x \in \setH^{M \times 1}
    \right \rbrace
  \\    
  \mathcal{RC}\left( \mat A \right) &:= 
    \left \lbrace
      \vec y \in \setH^{N \times 1}
      \colon
      \vecTr y = \vecTr x \mMulR \mat A, \forall \vec x \in \setH^{M \times 1}
    \right \rbrace  
\end{align*}

\subsection*{Left row null space ($\mathcal{LRN}$) and right row null space ($\mathcal{RRN}$)}
\begin{align*}
  \mathcal{LRN}\left( \mat A \right) &:= 
    \left \lbrace
      \vec x \in \setH^{M \times 1}
      \colon
      \mat A \mMulL \vec x = \vec 0
    \right \rbrace
  \\
  \mathcal{RRN}\left( \mat A \right) &:= 
    \left \lbrace
      \vec x \in \setH^{M \times 1}
      \colon
      \mat A \mMulR \vec x = \vec 0
    \right \rbrace  
\end{align*}

\subsection*{Left column null space ($\mathcal{LCN}$) and right column null space ($\mathcal{RCN}$)}
\begin{align*}
  \mathcal{LCN}\left( \mat A \right) &:= 
    \left \lbrace
      \vec x \in \setH^{N \times 1}
      \colon
      \vecTr x \mMulL \mat A  = \vecTr 0
    \right \rbrace
  \\
  \mathcal{RCN}\left( \mat A \right) &:= 
    \left \lbrace
      \vec x \in \setH^{N \times 1}
      \colon
      \vecTr x \mMulR \mat A  = \vecTr 0
    \right \rbrace  
\end{align*}

\subsection*{Relations between subspaces}

Similar to the complex case, the transpose operation relates some of the subspaces
to each other.

\begin{align*}
  \mathcal{LR}\left( \mat A \right) &= \mathcal{RC}\left( \matTr A \right) \\
  \mathcal{RR}\left( \mat A \right) &= \mathcal{LC}\left( \matTr A \right) \\
  \mathcal{LRN}\left( \mat A \right) &= \mathcal{RCN}\left( \matTr A \right) \\
  \mathcal{RRN}\left( \mat A \right) &= \mathcal{LCN}\left( \matTr A \right)  
\end{align*}
\section{The Matrix Inverse}

We now turn to investigate matrix inverses based on the proposed left and right
multiplication.

Let $\mat X$ be some matrix  $\mat X \in \setH^{M \times N}$ and let 
$\mat A \in \setH^{M \times M}$ be a square matrix. The left inverse $\mat A\mInvL$ 
satisfies the following condition.
\begin{align}
  \label{eq:LLinvCond}
  \mat A\mInvL \mMulL \left( \mat A \mMulL \mat X \right)
  = \mat X
\end{align}
Similarly, the right inverse $\mat A\mInvR$ satisfies this condition:
\begin{align}
\label{eq:RRinvCond}
  \mat A\mInvR \mMulR \left( \mat A \mMulR \mat X \right) 
  = \mat X    
\end{align}

By exploiting the associativity properties \gleichung{identA7} and \gleichung{identA8}
one can observe the following identities.
\begin{align*}
  \mat A\mInvL \mMulL \left( \mat A \mMulL \mat X \right) 
    &= \left(\mat A\mInvL \mMulL \mat A \right) \mMulL \mat X = \mat X \\
  \mat A\mInvR \mMulR \left( \mat A \mMulR \mat X \right) 
    &= \left(\mat A\mInvR \mMulR \mat A \right) \mMulR \mat X = \mat X    
\end{align*}
Hence, we have:
\begin{align}
  \label{eq:Alinv_ident}
  \mat A\mInvL \mMulL \mat A &= \mat I_M \\
  \label{eq:Arinv_ident}
  \mat A\mInvR \mMulR \mat A &= \mat I_M
\end{align}
In order to define the left inverse $\mat A\mInvL$ and right inverse 
$\mat A\mInvL$ we use the symplectic decomposition of $\mat A$.
\begin{align*}
  \mat A &:= \mat A_0 + \mat A_1 \mu_\perp, \quad \mat A_i\in \setC_\mu^{M \times M}
\end{align*}
Additionally, let us define the \emph{left adjoint matrix}\footnote{In \cite{Zhang97} and \cite{LeBih04}
this matrix is simply called \emph{adjoint matrix} since only left matrix multiplication has been considered.}
\begin{align*}
  \mat \chi_\mu\left\lbrace \mat A\right \rbrace &:= \begin{bmatrix} \mat A_0 & \mat A_1 \\ -\mat A_1\mStar & \mat A_0\mStar \end{bmatrix}
\end{align*}  
as well as the \emph{right adjoint matrix}
\begin{align*}
  \mat \chi'_\mu\left\lbrace \mat A\right \rbrace &:= \begin{bmatrix} \mat A_0 & -\mat A_1\mStar \\ \mat A_1 & \mat A_0\mStar \end{bmatrix} .
\end{align*}
These matrices render the direct complex-valued representation of a quaternion matrix. The left adjoint
matrix is connected to the left matrix multiplication whereas the right adjoint matrix inherently represents
the right matrix multiplication.
Therefore, the left and right inverses can be computed with complex arithmetics by noting the following two identities
(see also \cite{LeBih04}).
\begin{align*}
  \mat \chi_\mu\left\lbrace \mat A\mInvL \right \rbrace &= \mat \chi\mInv_\mu\left\lbrace \mat A \right \rbrace 
  \\
  \mat \chi'_\mu\left\lbrace \mat A\mInvR \right \rbrace &= \mat {\chi'}\mInv_\mu\left\lbrace \mat A \right \rbrace
\end{align*}  

Note that if $\mat A = \mat A_0$ and therefore $\mat A \in \setC_\mu^{M \times M}$ both adjoint matrices
are equal. That is, for fields which are isomorphic to the compelx numbers both, $\mat A\mInvL$ and $\mat A\mInvR$, 
reduce to the complex inverse $\mat A\mInv$.


Without giving the proof we state that the right and left inverse defined above
remain the same when the order of the matrix products is changed.
\begin{align*}
  \left( \mat X \mMulL \mat A \right) \mMulL \mat A\mInvL
  &= \mat X
  \\
  \left( \mat X \mMulR \mat A \right) \mMulR \mat A\mInvR
  &= \mat X    
\end{align*}
Hence, in addition to \gleichung{Alinv_ident} and \gleichung{Arinv_ident} we have
the following identities:
\begin{align*}
  \mat A \mMulL \mat A\mInvL &= \mat I_M \\
  \mat A \mMulR \mat A\mInvR &= \mat I_M
\end{align*}

We complete the investigation of the quaternion matrix inverses by pointing how
they are connected.
\begin{align}
	\label{eq:inv_relation}
	\left(\mat A\mInvL\right)\mTr	
	=
	\left(\matTr A\right)\mInvR
	\quad
	\Leftrightarrow
	\quad
	\left(\mat A\mInvR\right)\mTr
	=
	\left(\matTr A\right)\mInvL
\end{align}
This becomes clear by applying \gleichung{qmatTransp}.
\begin{align*}
	\mat I_M 
		&= \left( \mat A\mInvL \mMulL \mat A \right)\mTr
		\\
		& 
		=
		\left(\mat A\mInvL\right)\mTr \mMulR \matTr A
\end{align*}
Hence, $\matTr A$ must be the right inverse of $\left(\mat A\mInvL\right)\mTr$ and vice versa.
Moreover, if all entries of $\mat A$ are located in the same set $\setC_\mu^{M \times M}$,
\gleichung{inv_relation} reduces to the well known relation 
$\left(\mat A\mInv\right)\mTr = \left(\mat A\mTr\right)\mInv$.

\section{The Kronecker Product And The Khatri-Rao Product}

\subsection{The Left and Right Kronecker Product}

In addition to the matrix multiplication, it is also desirable to examine Kronecker
products of two quaternion matrices $\mat A \in \setH^{M_\mathrm A \times N_\mathrm A}$ and 
$\mat B \in \setH^{M_\mathrm B \times N_\mathrm B}$. 
In this case, one must also distinguish between the left Kronecker product 
\begin{align*}
  \mat A \mKronL \mat B
    := \begin{bmatrix}
	\left[ \mat A \right]_{1,1} \cdot \mat B &
	\dots &
	\left[ \mat A \right]_{1,N_\mathrm A} \cdot \mat B
	\\
	\vdots & \ddots & \vdots
	\\
	\left[ \mat A \right]_{M_\mathrm A,1} \cdot \mat B &
	\dots &
	\left[ \mat A \right]_{M_\mathrm A,N_\mathrm A} \cdot \mat B
      \end{bmatrix}
\end{align*}
and the right Kronecker product
\begin{align*}
  \mat A \mKronR \mat B
    := \begin{bmatrix}
	\mat B \cdot \left[ \mat A \right]_{1,1} &
	\dots &
	\mat B \cdot \left[ \mat A \right]_{1,N_\mathrm A}
	\\
	\vdots & \ddots & \vdots
	\\
	\mat B \cdot \left[ \mat A \right]_{M_\mathrm A,1}  &
	\dots &
	\mat B \cdot \left[ \mat A \right]_{M_\mathrm A,N_\mathrm A}
      \end{bmatrix} .
\end{align*}
In contrast to the matrix product, transposing the Kronecker product 
does not change its type. 
\begin{align*}
\left(\mat A \mKronL \mat B\right)\mTr 
  = \matTr A \mKronL \matTr B \\
\left(\mat A \mKronR \mat B\right)\mTr 
  = \matTr A \mKronR \matTr B  
\end{align*}

\subsection{Kronecker Product and Vectorization}

It is possible to reformulate the vectorization of a product of three matrices 
$\mat A \in \setH^{M \times K}$, $\mat B \in \setH^{K \times L}$, 
and $\mat C \in \setH^{L \times N}$ in the quaternion domain.
\begin{align*}
  \mVec{ \mat A \mMulL \left[ \mat B \mMulR \mat C \right] } &= 
    \left(\matTr C \mKronR \mat A \right) \mMulL \mVec{ \mat B }
  \\
  \mVec{ \mat A \mMulR \left[ \mat B \mMulL \mat C \right] } &= 
    \left(\matTr C \mKronL \mat A \right) \mMulR \mVec{ \mat B }
  \\
  \mVec{ \left[ \mat A \mMulL \mat B \right] \mMulR \mat C  } &= 
    \left(\matTr C \mKronL \mat A \right) \mMulL \mVec{ \mat B }  
  \\
  \mVec{ \left[ \mat A \mMulR \mat B \right] \mMulL \mat C  } &= 
    \left(\matTr C \mKronR \mat A \right) \mMulR \mVec{ \mat B }  
\end{align*}
However, there is no such expression for
$\mVec{ \mat A \mMulL \mat B \mMulL \mat C  }$ and
$\mVec{ \mat A \mMulR \mat B \mMulR \mat C  }$ using the operators
proposed in this work.
The reason is that this would involve multiplying the components of 
$\mat B$ between $\mat A$ and $\mat C$.

\subsection{The Khatri-Rao Product}

The last pair of operators we introduce is the left  Khatri-Rao product
\begin{align*}
  \mat C \mKhatriL \mat D :=
    \begin{bmatrix} 
      \vec c_1 \mKronL \vec d_1 &
      \dots &
      \vec c_N \mKronL \vec d_N 
    \end{bmatrix}
\end{align*}
as well as the right Kathri-Rao product
\begin{align*}
  \mat C \mKhatriR \mat D :=
    \begin{bmatrix} 
      \vec c_1 \mKronR \vec d_1 &
      \dots &
      \vec c_N \mKronR \vec d_N 
    \end{bmatrix}
\end{align*}
of two marices 
$\mat C = \begin{bmatrix} \vec c_1 & \dots & \vec c_N \end{bmatrix}$
and
$\mat D = \begin{bmatrix} \vec d_1 & \dots & \vec d_N \end{bmatrix}$.

\subsection{The Khatri-Rao Product and Vectorization}

Let $\mat A \in \setH^{M \times K}$ and $\mat C \in \setH^{K \times N}$
be two matrices. Additionally, let 
$
\mat B = \mDiag{\vec b}
$ 
be a diagonal matrix having the entries of a vector 
$\mat b \in \setH^{K \times 1}$
on its main diagonal.
In this case, the following identities hold.
\begin{align*}
  \mVec{ \mat A \mMulL \left[ \mat B \mMulR \mat C \right] } &= 
    \left(\matTr C \mKhatriR \mat A \right) \mMulL \vec b
  \\
  \mVec{ \mat A \mMulR \left[ \mat B \mMulL \mat C \right] } &= 
    \left(\matTr C \mKhatriL \mat A \right) \mMulR \vec b
  \\
  \mVec{ \left[ \mat A \mMulL \mat B \right] \mMulR \mat C  } &= 
    \left(\matTr C \mKhatriL \mat A \right) \mMulL \vec b
  \\
  \mVec{ \left[ \mat A \mMulR \mat B \right] \mMulL \mat C  } &= 
    \left(\matTr C \mKhatriR \mat A \right) \mMulR \vec b
\end{align*}
However, there is no such expression for
$\mVec{ \mat A \mMulL \mat B \mMulL \mat C  }$ and
$\mVec{ \mat A \mMulR \mat B \mMulR \mat C  }$.
\section{Examples}

In this section we will give some examples showing the
benefit of using the proposed ways of notation.

\subsection{Systems of Complex Widely Linear Equations}

Let us consider a widely linear system of equations in the complex domain
with 
$\mat A, \mat B \in \setC^{M \times M}$,
$\vec X \in \setC^{M \times P}$, and
$\vec C \in \setC^{M \times P}$,
\begin{align}
  \label{eq:ex1_cmplx}
  \mat A \mat X + \mat B \mat X\mStar = \mat C
\end{align}

The complex conjugate prevents factoring out  $\mat X$.
A possible solution would be to look at the real part and imaginary part
of the above equation and to compute the result in the real domain.

Nevertheless, we may also consider \gleichung{ex1_cmplx} as a quaternion
equation with 
$\mat A, \mat B \in \setC_\mu^{M \times M}$,
$\mat X \in \setC_\mu^{M \times P}$, and
$\mat C \in \setC_\mu^{M \times P}$,
where $\mu \in \setHpu$ is an arbitrary pure unit quaternion (PUQ). Then, let $\mu_\perp \in \setHpu$
be another PUQ that is orthogonal to $\mu$. It can readily be verified
that $\mat X\mStar = - \mu_\perp \mat X \mu_\perp$. Therefore we may 
rewrite \gleichung{ex1_cmplx} as follows.
\begin{align}
  \label{eq:ex1_quat1}
  \mat A \mMulL \mat X - \mat B \mMulL (\mu_\perp \mat X \mu_\perp) = \mat X
\end{align}
Using the left matrix product \gleichung{ex1_quat1} becomes:
\begin{align*}
  \begin{bmatrix}
    \mat A & -\mat B \mu_\perp
  \end{bmatrix}
  \mMulL
  \begin{bmatrix}
    \mat X \\
    \mat X \mu_\perp
  \end{bmatrix}
  =
  \mat C
\end{align*}
Next, by using the right matrix product we may factor out $\mat X$.
\begin{align*}
  \begin{bmatrix}
    \mat A & -\mat B \mu_\perp
  \end{bmatrix}
  \mMulL
  \left(
  \begin{bmatrix}
    \mat I_M \\
    \mu_\perp \mat I_M
  \end{bmatrix}
  \mMulR
  \mat X
  \right)
  =
  \mat C
\end{align*}
Moreover, let us define the matrices $\mat F_1$ and $\mat G$.
\begin{align*}
  \mat F_1 &:=
  \begin{bmatrix}
    \mat A & -\mat B \mu_\perp
  \end{bmatrix}
  &
  \mat G &:=
  \begin{bmatrix}
    \mat I_M \\
    \mu_\perp \mat I_M
  \end{bmatrix}
\end{align*}
Using these matrices we arrive at a compact expression.
\begin{align}
  \label{eq:part1_widely_lin}
  \mat F_1 \mMulL \left( \mat G \mMulR \mat X \right) = \mat C
\end{align}

The same derivation can be done for \gleichung{ex1_cmplx} having
taken the conjugate on both sides of the equation.
\begin{align*}
  \mat B\mStar \mat X + \mat A\mStar \mat X\mStar = \mat C\mStar
\end{align*}
The result is similar to \gleichung{part1_widely_lin}
\begin{align}
  \label{eq:part2_widely_lin}
  \mat F_2 \mMulL \left( \mat G \mMulR \mat X \right) = \mat C\mStar,
\end{align}
with
\begin{align*}
  \mat F_2 &:=
  \begin{bmatrix}
    \mat B\mStar & -\mat A\mStar \mu_\perp
  \end{bmatrix} .
\end{align*}

Let us put it all together by defining the matrices $\mat F$ and $\mat C_\mathrm a$. 
\begin{align*}
  \mat F &:=
  \begin{bmatrix}
    \mat F_1 \\ \mat F_2
  \end{bmatrix} 
  &
  \mat C_\mathrm a &:=
  \begin{bmatrix}
    \mat C \\ \mat C\mStar
  \end{bmatrix}   
\end{align*}
By noting that $\mat G$ has orthogonal columns
\begin{align*}
  \matHerm G \mMulL \mat G  
  =
  \matHerm G \mMulR \mat G  
  = 2 \mat I_M,
\end{align*}
the following solution to \gleichung{ex1_cmplx} is obtained
\begin{align*}
  \mat X = 0.5\ \matHerm G \mMulR \left(\mat F\mInvL \mMulL \mat C_\mathrm a \right)
\end{align*}

\subsection{The Eigendecomposition}

As mentioned in \abschnitt{introduction} there has already been some research concerning
the eigenvalues of a quaternion matrix $\mat A \in \setH^{M \times M}$. We will now review
the definition of quaternionic eigenvectors and eigenvalues. Based on this, we examine how 
the proposed notation renders useful when defining Eigenvalue decompositions.

\subsubsection*{The Left Eigendecomposition}

If for a scalar $\lambda_\mathrm L \in \setH$ and a vector $\vec q_\mathrm L \in \setH^{M \times 1}$
the following condition holds, $\lambda_\mathrm L \in \setH$ is called a \emph{left eigenvalue} and 
$\vec q_\mathrm L$ is said to be a {left eigenvector} of $\mat A$.
\begin{align*}
    \mat A \vec q_\mathrm L = \lambda_\mathrm L \vec q_\mathrm L
\end{align*}
Hence, we call 
\begin{align*}
  \mat A = \left( \mat Q_\mathrm L \mMulR \mat \Lambda_\mathrm L \right) \mMulL \mat Q_\mathrm L\mInvL
\end{align*}
the \emph{left eigendecomposition} of $\mat A$,
where the matrix 
$\mat Q_\mathrm L := \begin{bmatrix}\vec q_{\mathrm L, 1} & \dots & \vec q_{\mathrm L, M}\end{bmatrix}$
stores the left eigenvectors corresponding to the eigenvalues
stored in the diagonal matrix 
$\mat \Lambda_\mathrm L := \mDiag{\lambda_{\mathrm L, 1}, \ldots, \lambda_{\mathrm L, M}}$.

\subsubsection*{The Right Eigendecomposition}

If for a scalar $\lambda_\mathrm R \in \setH$ and a vector $\vec q_\mathrm R \in \setH^{M \times 1}$
the following condition holds, $\lambda_\mathrm R \in \setH$ is called a \emph{right eigenvalue} and 
$\vec q_\mathrm R$ is said to be a {right eigenvector} of $\mat A$.
\begin{align}
  \label{eq:right_EV}
    \mat A \vec q_\mathrm R = \vec q_\mathrm R \lambda_\mathrm R
\end{align}
However, in general there might exist infinitely many right eigenvalues. This can be seen by multiplying
\gleichung{right_EV} from the right with a non-zero quaternion $s \in \setH$.
\begin{align*}
    \mat A \vec q_\mathrm R s &= \vec q_\mathrm R \lambda_\mathrm R s  \\
    \mat A \underbrace{\vec q_\mathrm R s}_{\vec q'_\mathrm R}
      &= \underbrace{\vec q_\mathrm R s}_{\vec q'_\mathrm R} 
          \underbrace{s^{-1} \lambda_\mathrm R s}_{\lambda'_\mathrm R}  \\
  \mat A \vec q'_\mathrm R  &= \vec q'_\mathrm R \lambda'_\mathrm R           
\end{align*}
Hence, if $\lambda_\mathrm R$ is a right eigenvalue of $\mat A$, then all similar quaternions  
$\lambda'_\mathrm R \in \mathcal S\left( \lambda_\mathrm R\right)$ are right eigenvalues as well
(see \cite[p. 36]{Zhang97}).

To define a \emph{right eigendecomposition}, let  
$\mat \Lambda_\mathrm R := \mDiag{\lambda_{\mathrm R, 1}, \ldots, \lambda_{\mathrm R, M}}$
be a diagonal matrix of mutually non-similar right eigenvalues. Moreover, let 
$\mat Q_\mathrm R := \begin{bmatrix}\vec q_{\mathrm R, 1} & \dots & \vec q_{\mathrm R, M}\end{bmatrix}$
be the the matrix of corresponding right eigenvectors. The resulting eigendecomposition obtains the 
following form.
\begin{align*}
  \mat A =  \mat Q_\mathrm R \mMulL \mat \Lambda_\mathrm R  \mMulL \mat Q_\mathrm R\mInvL
\end{align*}

Note that this decomposition is not unique. However, according to \cite[theorem 5.4]{Zhang97}
any matrix $\mat A \in \setH^{M \times M}$ has exactly $M$ complex right eigenvalues.
Additionally, these eigenvalues all have non-negative imaginary parts. These eigenvalues are
called \emph{standard eigenvalues} and may be used to resolve the uniqueness problem.
\subsection{The Quaternion Discrete Fourier Transform}

The Quaternion Discrete Fourier Transform (QDFT) has already been used in
applications such as image processing (see \cite{Ell07}). However, it
still lacks a convenient notation.

Let us first have a look at the hitherto QDFT notation of a matrix $\mat A \in \setH^{M \times N}$.
To this end, let $f_{1, m, u}^{(\mu_1)}$ and $f_{2, n, v}^{(\mu_2)}$ denote the Fourier basis function
of the QDFT, where $\mu_i \in \setHpu$.
\begin{align}
  \label{eq:fourier_basis}
  f_{1, m, u}^{(\mu_1)} &:= \tfrac{1}{\sqrt{M}} \exp \left(-\mu_1 2 \pi u \tfrac{m}{M} \right) \\
  f_{2, n, v}^{(\mu_2)} &:= \tfrac{1}{\sqrt{N}} \exp \left(-\mu_2 2 \pi v \tfrac{n}{N} \right)
\end{align}

The two-side, left-side and right-side DQFT arise by multiplying the basis functions from
different directions with respect to the entries of $\mat A$ (see also \cite{Pei01}). 
\begin{align}
  \label{eq:dqft1_comp}
  \left[ \mathcal{F}_{\mu_1, \mu_2}^{(1)} \left \lbrace \mat A \right \rbrace \right]_{u, v}
  &=
    \sum_{m=0}^{M-1} \sum_{n=0}^{N-1} 
    f_{1, m, u}^{(\mu_1)}
    \cdot 
    \left[\mat A \right]_{m, n} 
    \cdot
    f_{2, n, v}^{(\mu_2)}
\\
  \label{eq:dqft2_comp}
  \left[ \mathcal{F}_{\mu_1, \mu_2}^{(2)} \left \lbrace \mat A \right \rbrace \right]_{u, v}
  &=
    \sum_{m=0}^{M-1} \sum_{n=0}^{N-1} 
    f_{1, m, u}^{(\mu_1)}
    f_{2, n, v}^{(\mu_2)} 
    \cdot \left[\mat A \right]_{m, n} 
\\
  \label{eq:dqft3_comp}
  \left[ \mathcal{F}_{\mu_1, \mu_2}^{(3)} \left \lbrace \mat A \right \rbrace \right]_{u, v}
  &=
    \sum_{m=0}^{M-1} \sum_{n=0}^{N-1}  
    \left[\mat A \right]_{m, n} 
    \cdot
    f_{1, m, u}^{(\mu_1)}
    f_{2, n, v}^{(\mu_2)} 
\end{align}
Notice that for each DQFT it would be possible to swap the order of the Fourier basis functions
However, we will not consider these DQFTs in the following.

The authors of \cite{Ell07} presented the special case where $\mu_1 = \mu_2$ for the left-side and right-side DQFT.
Additionally, the authors did not take advantage of a left matrix multiplication and a right matrix multiplication.
Hence, we present a new notation of the DQFT using the notation proposed in \abschnitt{left_right_mult}.
To that end, we define the quaternion Fourier matrices $\mat F_1^{(\mu_1)} \in \setC_{\mu_1}^{M \times M}$ and 
$\mat F_2^{(\mu_2)} \in \setC_{\mu_2}^{N \times N}$  based on \gleichung{fourier_basis} as follows.
\begin{align*}
 \mat F_1^{(\mu_1)} &:= \tfrac{1}{\sqrt{M}} \exp \left(-\mu_1 2 \pi \vec m \vecTr m M^{-1} \right) \\
 \mat F_2^{(\mu_2)} &:= \tfrac{1}{\sqrt{N}} \exp \left(-\mu_2 2 \pi \vec n \vecTr n N^{-1} \right)
\end{align*}
The index vectors  $\vec m$ and $\vec n$ are defined as
$\vec m = \left[ \begin{smallmatrix} 0 & \dots & M-1\end{smallmatrix} \right]\mTr$.
$\vec n = \left[ \begin{smallmatrix} 0 & \dots & N-1\end{smallmatrix} \right]\mTr$.
Now, the DQFTs given in \gleichung{dqft1_comp} -- \gleichung{dqft3_comp} of a matrix $\mat A \in \setH^{M \times N}$ 
may conveniently be written as follows.
\begin{align*}
  \mathcal{F}_{\mu_1, \mu_2}^{(1)} \left \lbrace \mat A \right \rbrace
    &= \mat F_1^{(\mu_1)} \mMulL \mat A \mMulL \mat F_2^{(\mu_2)}
    \\
  \mathcal{F}_{\mu_1, \mu_2}^{(2)} \left \lbrace \mat A \right \rbrace
    &= \mat F_1^{(\mu_1)} \mMulL \left( \mat A \mMulR \mat F_2^{(\mu_2)} \right)
    \\
  \mathcal{F}_{\mu_1, \mu_2}^{(3)} \left \lbrace \mat A \right \rbrace
    &= \left(\mat F_1^{(\mu_1)} \mMulR \mat A \right) \mMulL \mat F_2^{(\mu_2)}
\end{align*}
Moreover, the inverse discrete quaternion Fourier transform (IDQFT) of a matrix  $\mat A \in \setH^{M \times N}$
obtains the following form.
\begin{align*}
  \mathcal{\tilde F}_{\mu_1, \mu_2}^{(1)} \left \lbrace \mat A \right \rbrace 
    &= 	{\mat F_1^{(\mu_1)}}\mHerm 
	\mMulL 
	\mat A
	\mMulL 
	{\mat F_2^{(\mu_2)} }\mHerm
    \\
  \mathcal{\tilde F}_{\mu_1, \mu_2}^{(2)} \left \lbrace \mat A \right \rbrace 
    &= \left(
	  {\mat F_1^{(\mu_1)}}\mHerm
	  \mMulL
	  \mat A
      \right) 
      \mMulR
      {\mat F_2^{(\mu_2)}}\mHerm
    \\
  \mathcal{\tilde F}_{\mu_1, \mu_2}^{(3)} \left \lbrace \mat A \right \rbrace
    &= 	{\mat F_1^{(\mu_1)}}\mHerm
	\mMulR 
	\left( 
	  \mat A
	  \mMulL
	  {\mat F_2^{(\mu_2)}}\mHerm
	\right)
\end{align*}
In general, the matrix product defined by the left-side (I)DQFT and the right-side (I)DQFT 
is not associative. This changes when $\mu_1$ and  $\mu_2$ are chosen to be equal. In 
this case, the parenthesis can be set arbitrarily.

\section{Conclusions}

Using linear algebra over a skew-field such as the quaternions  often lacks a 
convenient notation. To this end, we propose the left and right matrix product. 
These products inherently cover the problem of defining the order in which the
matrix components are being multiplied.

Further on, it turns out that this notation is capable of defining the set of
quaternionic fundamental subspaces in a convenient manner. Additionally, the feasibility
of this notation is shown by applying it to the Eigenvalue decomposition as
well as to the Quaternion Discrete Fourier transform. In both cases a more
compact and simple form is achieved.

In complex linear algebra the Kronecker product as well as the Khatri-Rao 
product often prove beneficial. Hence, these products are shortly investigated
as well.

\section*{Acknowledgement}
The authors would like to thank Stephen J. Sangwine for his valuable comments which helped to improve this article.


\bibliographystyle{amsplain}
\bibliography{main}

\end{document}